\documentclass{article}
\usepackage{stywhispers,amsmath,epsfig}

\usepackage{amsmath}
\usepackage{amsthm}
\usepackage{amsfonts}
\usepackage{amssymb}
\usepackage{verbatim}
\usepackage{multirow}
\usepackage{graphicx}
\usepackage{subfigure}
\usepackage{graphicx}
\usepackage{algorithmic}
\usepackage{color}
\usepackage{cite}
\usepackage[linesnumbered,lined,ruled]{algorithm2e}
\usepackage{tikz}
\usepackage{epstopdf}
\usepackage{color}
\usepackage{booktabs}
\usepackage{threeparttable}
\usepackage{rotating}
\usepackage{footmisc}

\newcommand{\cblue}[1]{{\color{black} #1}}

\newcommand{\R}{\mathbb{R}}
\newcommand{\bY}{\textbf{Y}}
\newcommand{\bX}{\textbf{X}}
\newcommand{\bM}{\textbf{M}}
\newcommand{\bN}{\textbf{N}}
\newcommand{\bS}{\textbf{S}}
\newcommand{\bA}{\textbf{A}}

\newcommand{\ba}{\textbf{a}}
\newcommand{\bm}{\textbf{m}}
\newcommand{\bx}{\textbf{x}}

\newcommand{\norm}[1]{\left\|#1\right\|}

\title{Robust collaborative nonnegative matrix factorization for hyperspectra unmixing (R-CoNMF)}
\name{Jun Li$^1$,
       Jos\'e M. Bioucas-Dias$^2$,
        Antonio~Plaza$^3$,
        and Lin Liu$^1$\thanks{\cblue{This work
		was supported by  the Portuguese Science and Technology Foundation under
		Projects UID/EEA/50008/2013, andProject PTDC/EEI-PRO/1470/2012.}}}
\address{$^1$School of Geography and Planning,  Sun Yat-sen University,  China\\
$^2$Instituto de Telecomunica\c{c}\~{o}es,  Instituto Superior T\'ecnico, Universidade de Lisboa,  Portugal\\
$^3$Hyperspectral Computing Laboratory, Escuela Polit\'{e}cnica,
University of Extremadura, Spain.}
\begin{document}
\ninept
\maketitle
\begin{abstract}
\cblue{The recently introduced collaborative nonnegative matrix factorization (CoNMF) algorithm  was conceived to  simultaneously  estimate  the number of endmembers, the mixing matrix, and the fractional abundances  from hyperspectral linear mixtures.  This paper introduces R-CoNMF, which is a robust version of CoNMF. The  robustness  has been added   by a)  including a volume regularizer  which penalizes the distance to a mixing matrix inferred by a pure pixel algorithm; and by b)  introducing a new proximal alternating optimization (PAO)  algorithm for which convergence to a critical point is guaranteed. Our experimental results indicate that R-CoNMF  provides  effective estimates  both when the number of endmembers are unknown and  when they are known.}

\end{abstract}
\begin{keywords}
Hyperspectral imaging, spectral unmixing, endmember extraction, collaborative nonnegative matrix factorization (CoNMF).
\end{keywords}
\section{Introduction}
\label{sec:intro}

Spectral unmixing is one of the most active topics for remotely sensed hyperspectral data exploitation \cite{manolakis01,keshava02,bioucas2012hyperspectral, kenmaSPM}. In general, most methods use  linear mixture model (LMM) due to its simplicity, where LMM assumes that each observed pixel in a hyperspectral image is linearly combined by a set of  pure  spectra or \textit{endmembers}. Within the LMM paradigm, three main types of unmixing algorithms can be identified: geometrical, statistical, and sparse regression-based \cite{bioucas2012hyperspectral}. Geometrical unmixing algorithms work under the assumption that the endmembers of a hyperspectral image are the vertices of a simplex of minimum volume enclosing the dataset (i.e., the set of hyperspectral vectors)  or of a simplex of  maximum volume  contained in the  convex hull of the dataset \cite{craig}, such as the minimum-volume enclosing simplex (MVES) \cite{specialissue-chan-ma}, minimum volume simplex analysis (MVSA) \cite{mvsaJournal}, the successive volume maximization (SVMAX) \cite{specialissue-chan-ma}, and etc. Among minimum volume algorithms, we highlight the minimum volume constrained nonnegative matrix factorization (MVC-NMF \cite{miao07}) which is also based on NMF principles. As their name suggests, statistical methods \cite{bioucas2012hyperspectral} are based on analyzing mixed pixels by means of statistical principles, such as Bayesian approaches \cite{bioucasdirichlet,eches}. Finally, sparse regression-based algorithms \cite{iordache2011sparse} are based on expressing each mixed pixel in a scene as a linear combination of a finite set of pure spectral signatures that are known \textit{a priori} and available in a library. Although each of these methods exhibits their own advantages and disadvantages, the fact is that geometrical approaches have been the ones most frequently used by the hyperspectral research community up to now \cite{bioucas2012hyperspectral}. This is mainly due to their reduced --although still quite high-- computational cost when compared with the other types of unmixing algorithms, as well as to the fact that they represent a straightforward interpretation of the LMM. In order to fully unmix a given hyperspectral image by means of a geometrical method, the majority of the state-of-the-art approaches are based on dividing the whole process into three concatenated steps \cite{bioucas2012hyperspectral}: 1) estimation of the number of endmembers; 2) identification of the spectral signatures of these endmembers; and 3) estimation of the endmember abundances in each pixel of the scene. In the last few years, several techniques have been developed for addressing each part of the chain, with particular emphasis on the identification of endmembers (with and without assuming the presence of pure spectral signatures in the input hyperspectral data \cite{bioucas2012hyperspectral,kenmaSPM}).

 \cblue{The recently introduced CoNMF algorithm  \cite{conmf}  addresses simultaneously the  three unmixing stages. In addition to a quadratic  data term, CoMNF} uses two regularizers:  a)  the $\ell_{2,1}$ mixed norm  \cite{R2} applied to the abundance matrix, which promotes sparsity among the \cblue{rows  of that matrix, and, therefore, it selects the active endmembers, and b) a volume regularizer  which simultaneously promotes minimum volume and pushes the endmembers with zero abundances  toward the mean value of the dataset. In spite of the good results obtained with CoNMF,  it has two   weaknesses: a)  there is no guarantee of convergence;  b)  in some cases,   we observe a negative joint effect of the  volume and  mixed norm  regularizers. In this paper, we  modify CoMNF aiming at  removing those weaknesses. The  former is removed by  introducing a PAO  solver to compute a critical point of the objective function. The latter is removed by modifying the volume regularizer aiming at robustness to noise and degenerated simplexes as  often happens with real data and  relatively high number of endmembers.}

The remainder of this paper is organized as follows. \cblue{Section \ref{sec:approach} introduces R-CoNMF.  Section \ref{sec:results} presents results  based on   synthetic hyperspectral data sets.} Comparisons with the MVC-NMF are also included. Finally, Section \ref{sec:conclusions} draws some conclusions and hints at plausible future research lines.

\section{The proposed approach}
\label{sec:approach}

Let ${\bf Y}\equiv[{\bf y}_1,\dots,{\bf y}_n]\in\mathbb{R}^{d\times n}$  be matrix representation of  a hyperspectral image with $n$ spectral vectors and $d$ spectral bands.  Under the LMM, we have \cite{bioucas2012hyperspectral}:
\begin{align}\label{eq:LMM}
   \bY  =&  \;\bM\bS +\bN\\
       & \text{s.t.:}\quad\bS\geq0, \;\;\;{\bf 1}^T_n\bS = {\bf 1}_p
       \nonumber
\end{align}
where $\bM\equiv  [\bm_1,\dots,\bm_p]\in\R^{d\times p}$ is a so-called mixing matrix containing $p$ endmembers, $\bm_i$ denotes the $i$-th endmember signature, ${\bS} \equiv [{\bf s}_1,\dots,{\bf s}_n]\in\R^{p\times n}$ is the abundance matrix containing the endmember fractions ${\bf s}_i$,  for   pixels $i=1,\dots,n$,  $\bS \geq 0$ is the abundance non-negativity constraint (to be understood in component-wise sense), \cblue{and ${\bf 1}^T_n\bS = {\bf 1}_p$ is the sum-to-one constraint (${\bf  1}_p$ stands for a column vector of ones with $p$ elements)}.
Finally, $\textbf{N}$ collects the errors that may affect the measurement process (e.g., noise).

\subsection{Estimation criterion}
\begin{figure}
  \centering
  \includegraphics[width=8cm]{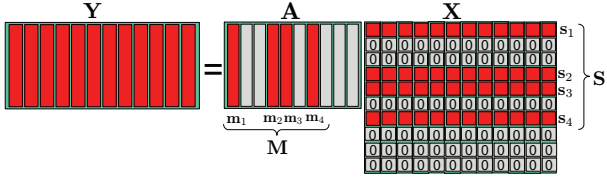}\\
  \caption{Illustration of the concept of  row  (collaborative) sparsity promoted $\ell_{2,1}$ norm under the linear mixing model. The abundance matrix $\bS$ is formed my the nonzero rows of $\bX$, while the mixing matrix $\bM$ is formed by the correspondent columns of $\bA$.}
  \label{fig:illustration_row_sparsity}
\end{figure}

In this work,  we address the estimation of $p$, {\em i.e}, the number of endmembers,  as well as the estimation of the mixing matrix $\bM$, and the abundance matrix $\bS$.  Although $p$ is not  known beforehand, we assume that we have access to an overestimate thereof. That is, we are given a number $q$ such that $q\geq p$. In this way, we account for a common  situation in which an overestimate of the number of endmembers is easy to compute, which is not usually the case  regarding the true  number of endmembers. \cblue{We tackle the estimation of $p$, $\bM$, and $\bS$ by seeking a solution for the following   NMF optimization:}
\begin{align}\label{eq:NMF_criterion}
     & \min_{\bf A,X}   \;\; (1/2) \norm{ \bY -\bA\bX}_F^2 + \alpha\norm{\bX}_{2,1}
           + \frac{\beta}{2}\norm{{\bf A}-{\bf P}}_F^2\\
          & \;\;\text{s.t.:}\quad\bX \in {\cal S}_ {q-1} \;\; {\bf   A} \in {\cal A}_{q-1} ,
          \nonumber
\end{align}
\cblue{where $\norm{\cdot}_F$ stands for the  Frobenius norm, $\bA \equiv [\ba_1,\dots,\ba_q]\in \R^{d\times q}$ and $\bX  \in \R^{q\times n}$ (${\bx}^i$ and ${\bx}_j$  will denote, respectively, the $i$-th row,  and the $i$-th column of $\bX$) are optimization variables, linked with the mixing matrix $\bf M$ and the abundance matrix $\bf S$, respectively. The relation between  $({\bf A},{\bf X})$ and $({\bf M},{\bf S})$ is illustrated in Fig. \ref{fig:illustration_row_sparsity}; the  term  $\norm{\bX}_{2,1}\equiv\sum_{i=1}^{q}\norm{\textbf{x}^i}_2 $ denotes the $\ell_{2,1}$ (see, {\em e.g.}, \cite{eldar2010average}) mixed norm of  matrix $\bf X$, the term
 $\norm{{\bf A}-{\bf P}}_F^2$ denotes the minimum volume regularizer, ${\bf P}\equiv[{\bf y}_{i_1},\dots,{\bf y}_{i_q}]$ is  a set of $q$  spectral   observed vectors
 inferred with  a pure-pixel algorithm, thus {\em close} to the extremes of the simplex,  $\alpha$ and $\beta$   are  regularization parameters, ${\cal  S}_ {q-1}$ is the collection of
 $q\times n$ matrices whose columns belong to the probability simplex of dimension $q-1$, and  $ {\cal A}_{q-1}$ is the collection of matrices of size $d\times q$ whose colums belong to the
 the affine set of dimension $q-1$ that best represents the data $\bf Y$ in the mean square error  sense.  The introduction of this constraint  removes the  shortcomings  associated violations to the sum-to-one constraints usually observed in real datasets.  }

\begin{figure}
  \centering
 \begin{tabular}{cc}
  \includegraphics[width=4cm]{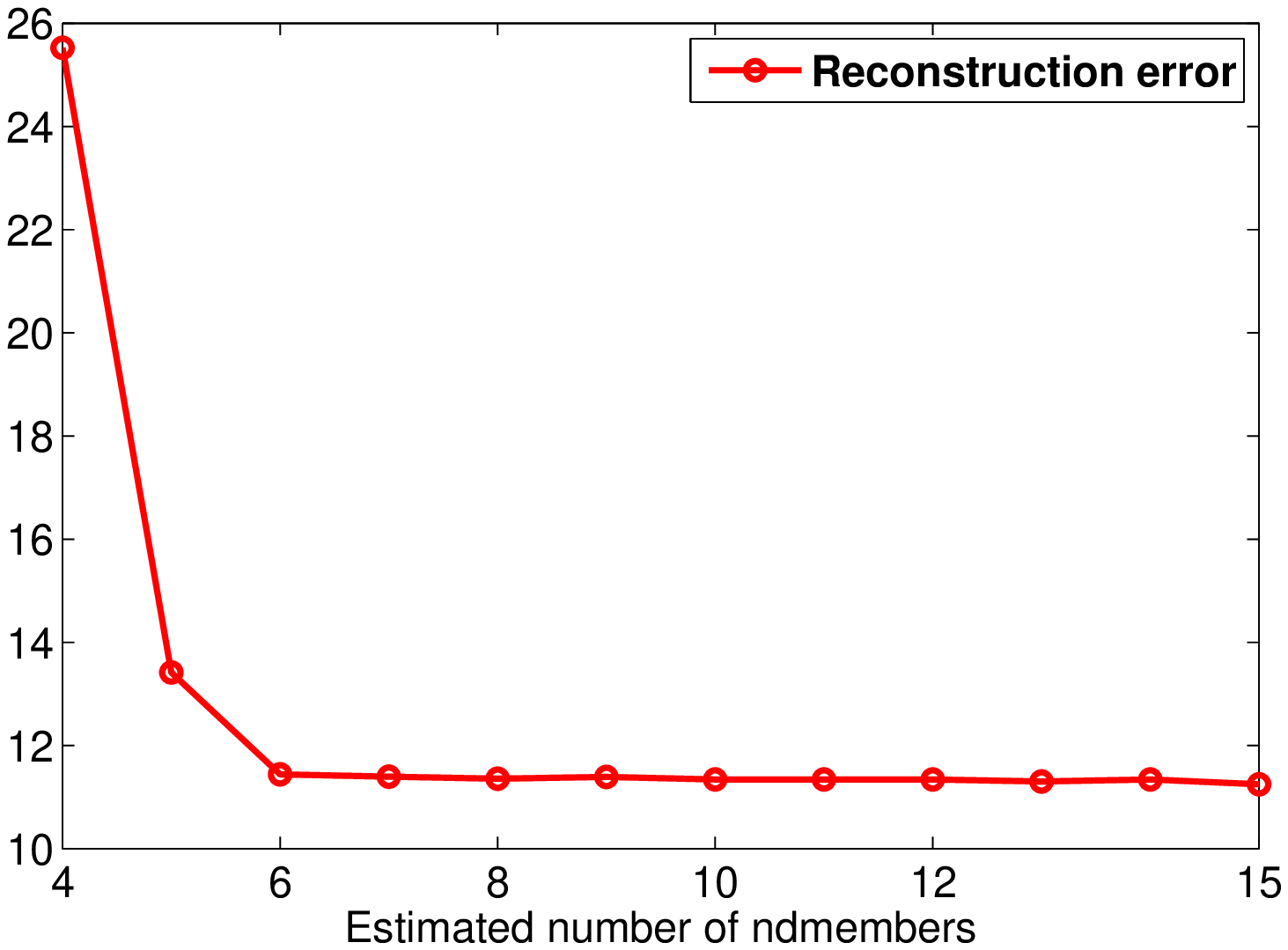} &
  \includegraphics[width=4cm]{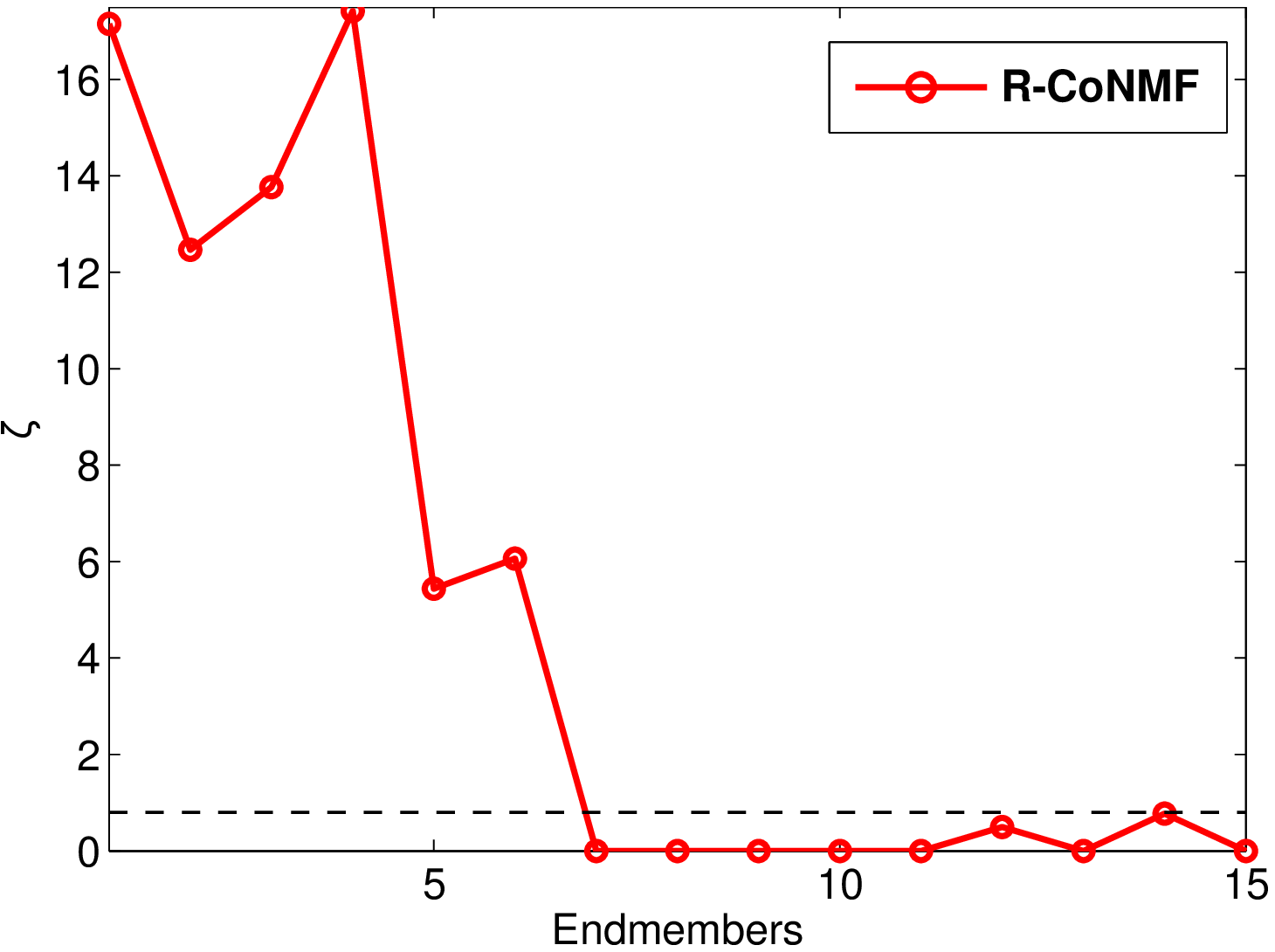} \\
(a) & (b) \\
  \end{tabular}
   \caption{(a) Reconstruction error as a function of the estimated number of endmembers. (b) Degree of sparseness $\zeta$ for $q=15$.}
  \label{fig:re}
\end{figure}

\cblue{
The objective function shown in  (\ref{eq:NMF_criterion}) has three terms: a  data fidelity term $\norm{\bY -\bA\bX}_F^2$, which promotes solutions  with low reconstruction error; an  $\ell_{2,1}$ mixed norm $\norm{\bX}_{2,1}$, which promotes row sparsity on $\bX$ \cite{fornasier2008recovery} (that is, it promotes solutions with complete rows $\bx^i $ set to zero); and a minimum volume term $\norm{{\bf A}-{\bf P}}_F^2$, which pushes the simplex defined by $\bf A$ towards the simplex defined by $\bf P$ which, given that  the it   is defined by observed vectors,
is inside the   true simplex.  This  volume regularizer   is a fundamental device of  R-CoNMF as it  largely reduces the sensitive of the solutions of \eqref{eq:NMF_criterion}
due  to  bad conditioning of the true mixing matrix $\bf M$ and to perturbations of the samples close to the  simplex facets.

\subsection{Optimization algorithm and convergence}
The nonconvex data fidelity  term $(1/2) \norm{\bY -\bA\bX}_F^2$ present in (\ref{eq:NMF_criterion})  makes the respective optimization hard.  Herein, we adopt the  following  PAO  iterative
algorithm:
\begin{align}
  \label{eq:opt_A}
      \bA_{(t+1)} &= \arg\min_{\bA \in {\cal A}_{q-1}} L(\bA, \bX_{(t)}) + \frac{\lambda}{2}\|\bA- \bA_{(t)}\|_F^2,\\
        \label{eq:opt_X}
      \bX_{(t+1)} & = \arg\min_{\bS \in {\cal S}_ {q-1}} L(\bA_{(t+1)}, \bX)+ \frac{\mu}{2}\|\bX- \bX_{(t)}\|_F^2,
\end{align}
where $L(\bA, \bX)$ is the objective function shown in \eqref{eq:NMF_criterion} and $\lambda,\mu > 0$ are two positive constants. We remark that  the above  procedure can be interpreted as a regularized version of a two block non-linear Gauss-Seidel method \cite{grippof1999globally}.  The POA algorithm \eqref{eq:opt_A}-\eqref{eq:opt_X} is an instance of the class considered in
\cite{Attouch:2010:PAM} for which the convergence to  a critical point   is proved in the Theorem 9 of the same paper. Optimization \eqref{eq:opt_A}  is a small size  quadratic problem, thus very light, and optimization \eqref{eq:opt_X} is a constrained $\ell_2-\ell_{2,1}$ problem which we solve effectively  with the CSUNSAL algorithm \cite{CosparseU}.

\subsection{Applying R-CoNMF }
R-CoNMF may be applied  either assuming that a)  the number of endmembers is know or b)  unknown. In the former case, we set $q=p$ and $\alpha$ to a very small value, thus removing the
$\ell_{2,1}$ regularizer from the objective function, although keeping the  constraint set ${\cal S}_ {q-1}$. In the  latter case, we  first apply R-CoNMF to infer the number of endmembers
and then we apply again R-CoNMF  as described in scenario a).

 Let us consider the scenario b). In this case, we run  R-CoNMF  for a fixed  $q > p$. Lets us define $\zeta(i)=\|{\bf x}^i\|_2$, for $i=1,\dots,q$, as a measure of the sparsity of abundances associated with the corresponding  endmembers ${\bf a}_i$, for  $i=1,\dots,q$. Ideally, we should  have $\zeta(i) = 0$ if ${\bf a}_i$ is not active.
  Due to the impact or noise and model errors,  we relax that criterion  as follows:  we  consider that an endmember ${\bf a}_i$  is active,  if  $\zeta(i) > \xi$, for a small  $\xi >0$. Fig. \ref{fig:re} shows the obtained $\zeta$ and the reconstruction error for a problem with $p=6$, $n=4000$, and  zero-mean Gaussian iid noise with SNR$=30$dB, where $\text{SNR}=10\log_{10}(\mathbb{E}[\|{\bf MS}\|_F^2/\mathbb{E}\|{\bf N}\|^2_F)$.  The applications of the above criterion  with any value of $\xi$ between 0.5 and 4 yields  $p=6$, which is the correct estimate of the number of endmembers. Notice that, this number could also be obtained from the plot of the reconstruction error. In  more complex scenarios with lower SNR and a larger number of endmembers, we may devise a strategy to combine both indicators.
  }

\section{Experimental results}
\label{sec:results}

In this section, we evaluate the proposed R-CoNMF method using synthetic hyperspectral data. The advantage of using synthetic scenes is that they offer a fully controlled analysis scenario in which the properties of our algorithm can be investigated precisely. The synthetic scenes have been generated using the LMM in (\ref{eq:LMM}). {The scenes comprise $n=4000$ pixels, and the spectral signatures used for their generation were randomly selected from the United States Geological Survey (USGS) digital spectral library\footnote{Available online: http://speclab.cr.usgs.gov/spectral-lib.html}. In order to ensure the {difference} among the endmembers used for simulation purposes,  the spectral angle distance (SAD) between any two spectral signatures is bigger than 10 degrees. Furthermore, let $p_{\text{mix}}$ be number of endmembers in one pixel. In real scenarios, it is possible to have a  large number of endmembers in a scene, for instance, $p\geq 10$. However, it is unlikely to have a large size of endmembers in one pixel. That is, in general, $p_{\text{mix}}$ is relatively small, say $p_{\text{mix}}\leq 5$. Based on this observation,   for the simulated data if $p\geq5$, we set $p_{\text{mix}}=5$. Otherwise, if $p<5$, we set $p_{\text{mix}}=p$. Finally,} to ensure that no pure pixels are present in the synthetic images, we discard all pixels with abundance fractions larger than 0.8, \emph{i.e.}, $\max({\bf s}_i) \leq 0.8$. 

\cblue{
In the case of $q=p$, let $\widehat{\bf M} = {\bf A}$ and  $\widehat{\bf S}= {\bf X}$ denote the  estimates of $\bf M$ and ${\bf S}$, respectively.} As performance indicators, we use {the relative reconstruction error $\text{RRE}=\|{\bf Y} - \widehat{\bf M}\widehat{\bf S}\|_F^2{/}\|{\bf Y}\|_F$,}  the SAD {(in degrees)}, and two error metrics focused on evaluating the quality of the estimated endmembers: $\|\widehat{\bf M}-{\bf M}\|_F$, and the quality of the estimated abundances: $\frac{1}{\sqrt{n \times p}}\|\widehat{\bf S}-{\bf S}\|_F$. 
It should be noted that, in all our experiments, we projected the data into a $q$-dimensional subspace with $q\geq p$.   \cblue{In Subsection \ref{exp1}, where the number of endmembers is assumed unknown, we set $\alpha = 10^{-5}$ and $\beta = 10^{-1}$. In Subsection \ref{exp2}, where the number of endmembers is assumed
known, we set $\beta = 10^{-5}$ and $\alpha$  was  hand tuned for optimal performance.}

\subsection{Experiment 1}
\label{exp1}

\begin{figure}
  \centering
 \begin{tabular}{cc}
   \includegraphics[width=4cm]{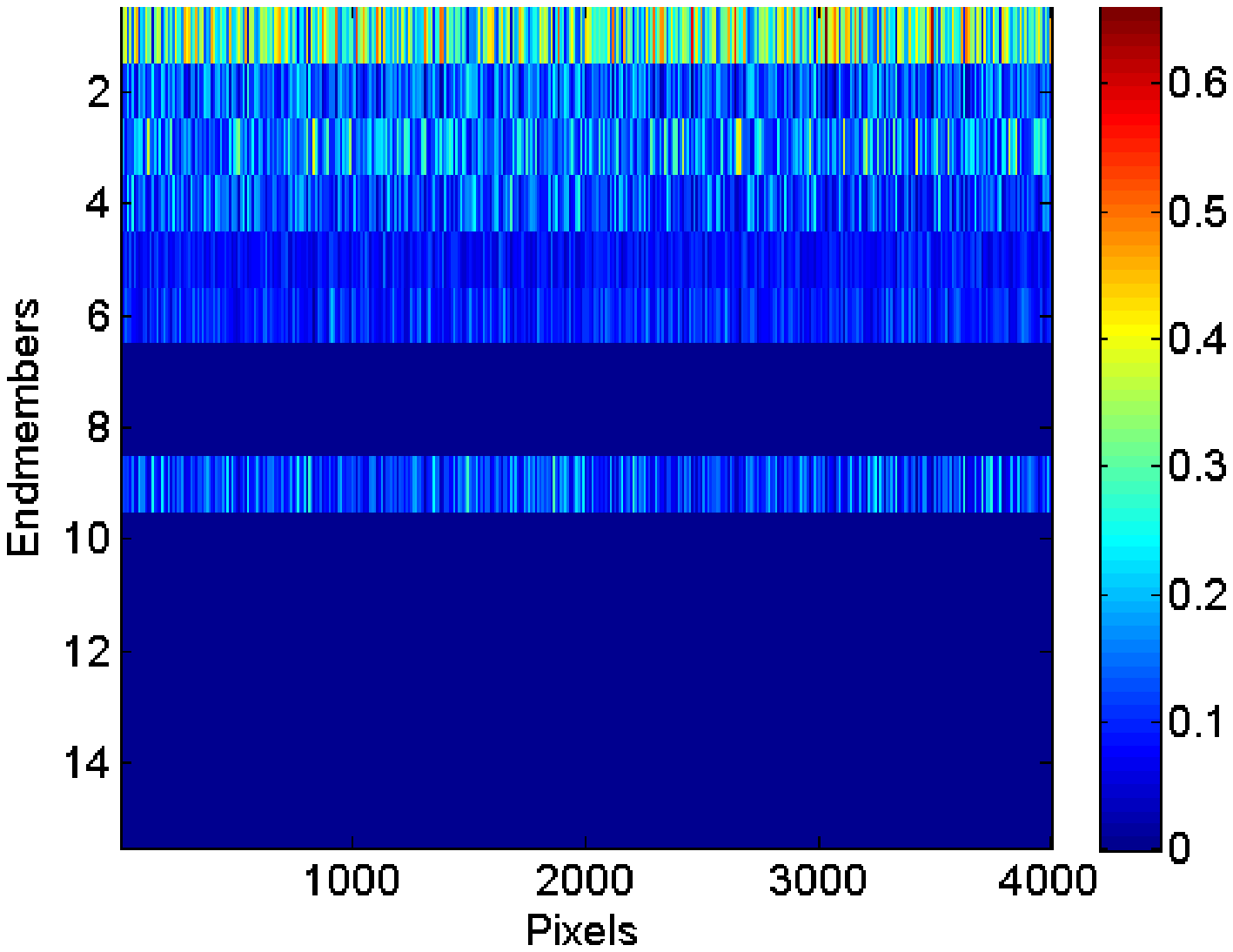} &
  \includegraphics[width=4cm]{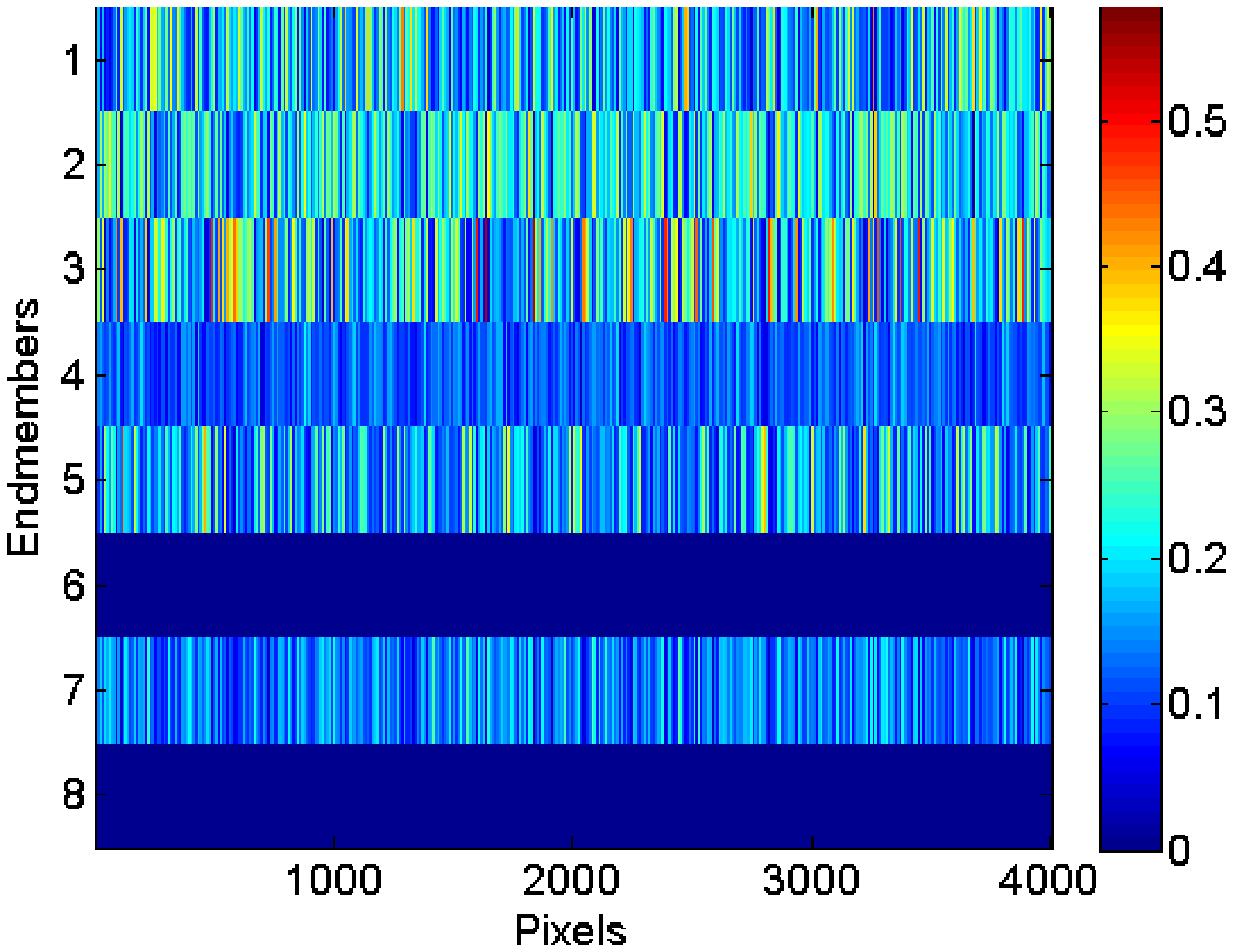} \\
  (a) $p=6$, $q=15$  &  (b) $p=6$, $q=8$\\

  \includegraphics[width=4cm]{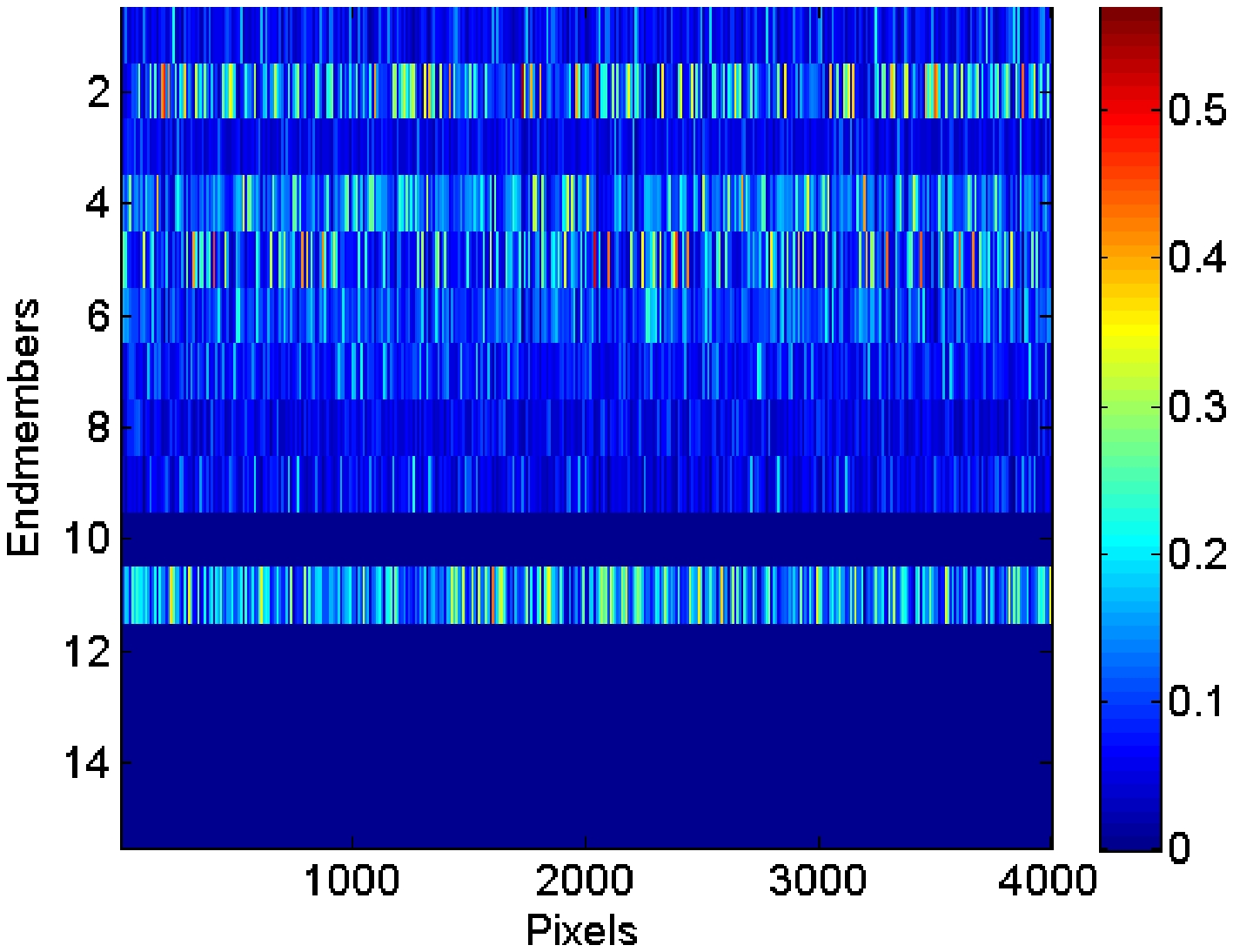} &
  \includegraphics[width=4cm]{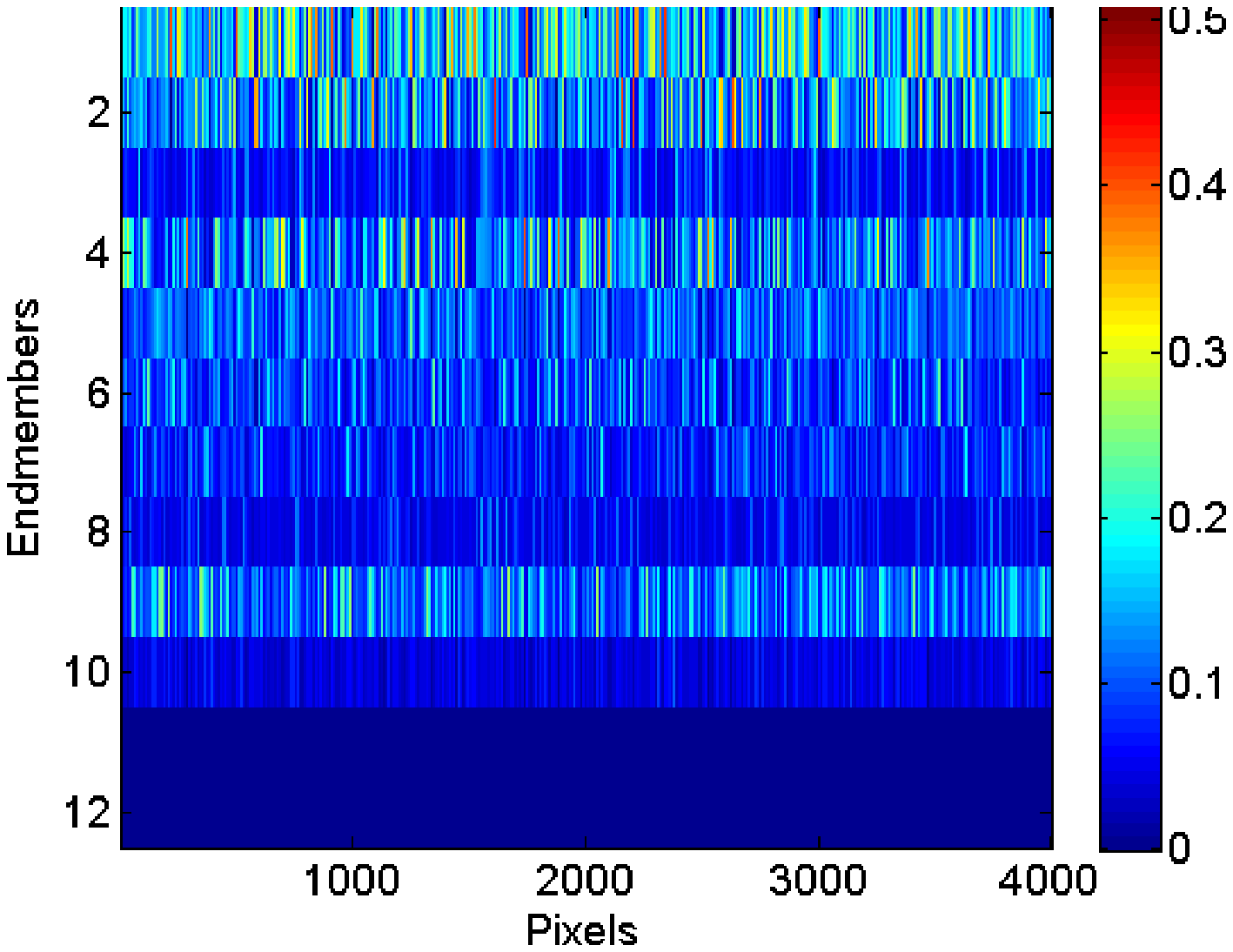} \\
  (c) $p=10$, $q=15$  &  (d) $p=10$, $q=12$\\

  \end{tabular}
   \caption{Estimated fractional abundance matrix $\widehat{\bf X}$ for problems with $n=4000$ pixels and SNR=30dB. }
  \label{fig:X}
\end{figure}

The first experiment aims at showing the good capability of R-CoNMF to provide the correct number of endmembers. Fig. \ref{fig:X} shows the obtained fractional abundance matrix for four different problems with $p=6$ and $p=10$ endmembers, respectively, by using different values of $q$. In all cases, we considered a scenario with SNR=30dB. Let $\widehat{q}$ be the number of endmembers estimated by R-CoNMF. In Figs. \ref{fig:X}(a) and \ref{fig:X}(b), where the number of endmembers is relative small (\emph{i.e.}, $p=6$), it is easy to detect the real number of endmembers, that is, $\widehat{p} = 6$. When the number of endmember increases, as shown in Figs. \ref{fig:X}(c) and \ref{fig:X}(d), the R-CoNMF can still provide a good estimate $\widehat{p} = 10$ even though the problem is much more difficult in this case.  In this experiment R-CoNMF provides effective estimates and works  according to our expectation. \cblue{However,  a  more detailed evaluation demands   extensive experiments with simulated and real data.}

\begin{figure}\scriptsize
\centering
\begin{tabular}{c}
\includegraphics[height=2.0in]{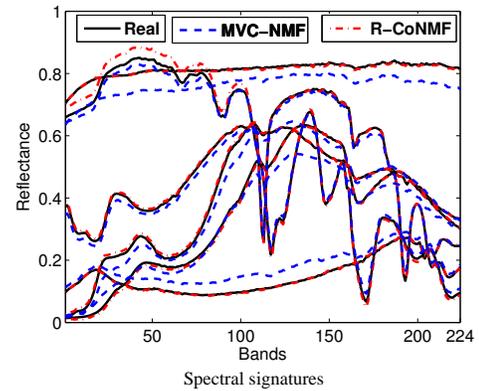}\\
 Spectral signatures \\
\end{tabular}
\caption{Spectral signatures of the endmembers extracted by R-CoNMF and MVC-NMF as compared to the reference signatures used for the simulation of a synthetic scene with $q=p=6$ and SNR=30dB.}
\label{p6pr6figurea}
\end{figure}

\begin{figure}\scriptsize
\centering
\begin{tabular}{cc}

\includegraphics[height=1.25in]{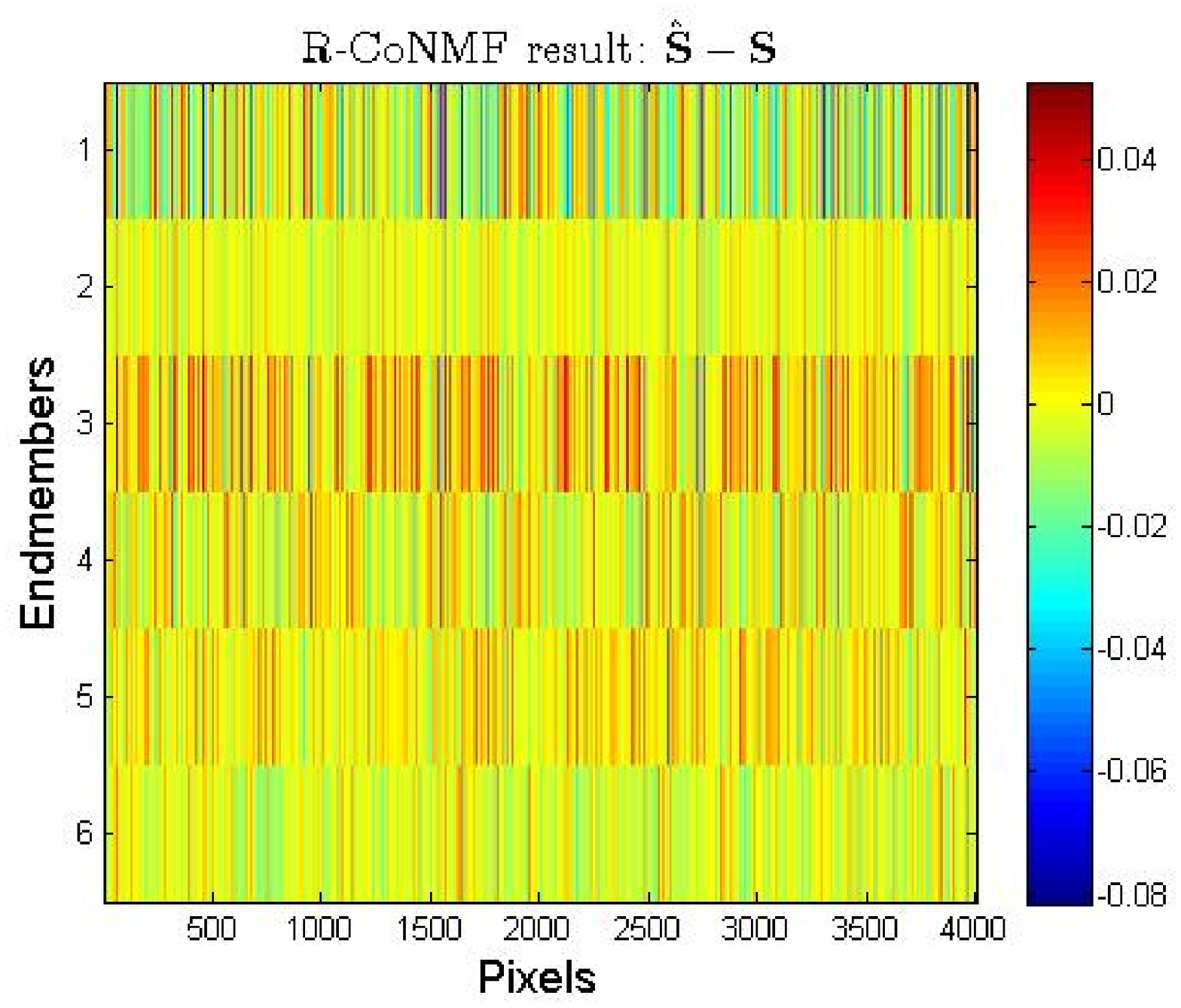} &
\includegraphics[height=1.25in]{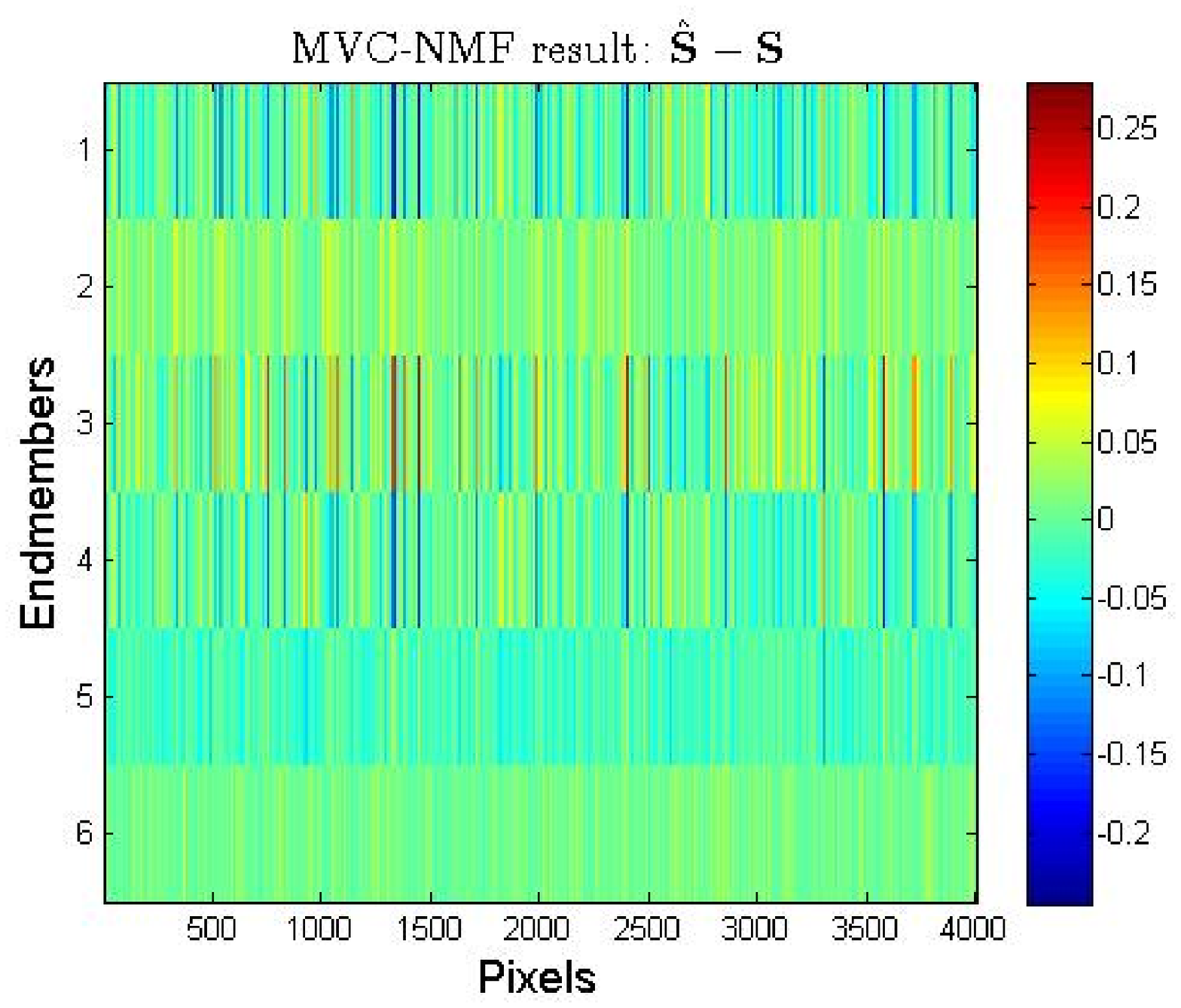}\\
(a) R-CoNMF  & (b) MVC-NMF\\
\end{tabular}
\caption{ (a) Difference  between real and estimated abundances for the R-CoNMF algorithm. (b) Difference  between the real and estimated abundances for the MVC-NMF algorithm.}
\label{p6pr6figureb}
\end{figure}

\begin{table*}\scriptsize
\renewcommand{\arraystretch}{1.3}
\caption{Evaluation of the performance of R-CoNMF and MVC-NMF in the unmixing of a synthetic hyperspectral data set, simulated with SNR=30dB, for different numbers of endmembers and $q = p$, where ``-'' means no results. }
\label{p4pr4table}
\begin{center}
\begin{tabular}{c|cccc|cccc}
\hline
\hline
\multirow{3}{*}{$q=p$} &   \multicolumn{4}{c|}{R-CoNMF}           &   \multicolumn{4}{c}{MVC-NMF} \\
   \cline{2-9}
 &  $\|\widehat{\bf M}-{\bf M}\|_F$   &  $\frac{1}{\sqrt{n \times p}}\|\widehat{\bf S}-{\bf S}\|_F$    & SAD  & RRE    
 &$\|\widehat{\bf M}-{\bf M}\|_F$   &  $\frac{1}{\sqrt{n \times p}}\|\widehat{\bf S}-{\bf S}\|_F$    & SAD  & RRE  \\ 
\hline
 4 & 0.17$\pm$0.09  &    0.01$\pm$3e-3 &    0.64$\pm$0.46 &   0.03$\pm$0.0 &
0.84$\pm$0.47 &   0.03$\pm$0.02   &  2.54$\pm$1.73 &   0.20$\pm$0.05 \\

6   &
 0.20$\pm$0.08 &   0.01$\pm$4e-3 &    0.56$\pm$0.21 &   0.03$\pm$1e-4 &
  1.44$\pm$1.17 &  0.04$\pm$0.02 &   3.48$\pm$2.57 &   0.23$\pm$0.04 \\

8  &  0.76$\pm$0.23 &    0.02$\pm$0.01 & 1.89$\pm$0.67 &    0.03$\pm$1e-3&
 36.50$\pm$154.23 &   0.08$\pm$0.05  &   11.97$\pm$19.48 &    0.25$\pm$0.07\\

10 &      1.28$\pm$0.50  &    0.03$\pm$0.01 &    2.77$\pm$1.73 &   0.03$\pm$1e-3 &
-  &	- & - & - \\

15   &   3.13$\pm$1.81 &    0.05$\pm$0.01 &    5.25$\pm$3.87 &    0.04$\pm$1e-3 &-  &	- & - & -  \\

\hline
\hline
%
%
%
%
%
\end{tabular}
\end{center}
\end{table*}

\subsection{Experiment 2}
\label{exp2}

In a second experiment we evaluate the performance of R-CoNMF in the case that the number of estimated endmembers coincides with the number of real endmembers, \emph{i.e.},  $q = p$. Here, we use the MVC-NMF algorithm in \cite{miao07} as a baseline for comparison with our method. Table \ref{p4pr4table} displays the results obtained by MVC-NMF and our proposed R-CoNMF algorithm for all considered performance discriminators for different values of $q = p = \{4,6,8,10,15\}$. In all cases, we considered an SNR of 30 dB and reported the results obtained from averaging the results of 30 Monte Carlo runs.

From the results reported on Table \ref{p4pr4table}, we can make the following observations. First and foremost, when there are only a few endmembers existing in the image (say, $q = p = 4$), both R-CoNMF and MVC-MNF obtained very good results. This is expected, since in this case it is relatively easy to solve the optimization problem. It is interesting to observe that, as the number of endmembers increases, the R-CoNMF  obtained very good performance (note the good performance obtained for the case $q = p = 10$). Even in a very difficult scenario such as $q = p = 15$, the solution provided by R-CoNMF is still useful. It should be noted that, in cases with a relatively high number of endmembers (i.e., $q = p \geq 8$), the MVC-NMF yields useless results. This is because, when the number of endmember increases, {most pixels are likely to fluctuate around the facets, a situation in which minimum volume-based algorithms are likely to fail \cite{bioucas2012hyperspectral}. } Even in this difficult scenario, in which the MVC-NMF could not provide feasible results (labeled as ``-" in Table \ref{p4pr4table}), the proposed R-CoNMF was able to provide a reasonable solution.  Based on this experiment, we can conclude that R-CoNMF is quite robust and has no strong constraints related with the quality of the analyzed data set.

{For illustrative purposes, Fig. \ref{p6pr6figurea} shows the signatures  estimated by R-CoNMF  and MVC-NMF. The  estimated spectral signatures by R-CoNMF are  similar to the real ones, while those estimated by MVC-NMF are slightly different. Similar observations can be made from the difference maps between the real and estimated abundance maps, as shown in Fig. \ref{p6pr6figureb}, where the difference of  R-CoNMF is much smaller  than that of the MVC-NMF.} The results in this subsection indicate that the R-CoNMF can perform very accurately in the case that the number of endmembers is known \textit{a priori}, i.e. $q = p$.

\section{Conclusions and future lines}
\label{sec:conclusions}

In this paper, we {propose} R-CoNMF, which is robust  version of the collaborative nonnegative matrix factorization (CoNMF) algorithm, which estimates  the number of endmembers, the  mixing matrix, and  the corresponding abundances. The proposed R-CoNMF algorithm fills a gap in the hyperspectral unmixing literature as it is one of the few algorithms that can accomplish the three main steps of the unmixing chain in fully automatic fashion, without the need to resort to external algorithms.  In the future, we will evaluate the algorithm using real hyperspectral data sets.

\bibliographystyle{IEEEbib} \scriptsize

\end{document}